\documentclass{amsart}

\usepackage[numbers]{natbib}
\usepackage{amsmath}
\usepackage[T1]{fontenc}
\usepackage{textcomp}
\usepackage{amssymb}
\usepackage{mathptmx}
\usepackage[numbers]{natbib}
\usepackage{bibunits}  
\usepackage{url}
\usepackage[bookmarks=true,pdfstartview={FitBH},colorlinks={true},
  pdfview={FitBH}]{hyperref}  

\theoremstyle{plain}
\newtheorem{thm}{Theorem}[section]
\newtheorem{claim}[thm]{Claim}
\newtheorem{lemma}[thm]{Lemma}
\newtheorem{corollary}[thm]{Corollary}
\newtheorem{proposition}[thm]{Proposition}

\theoremstyle{definition}
\newtheorem{definition}[thm]{Definition}

\theoremstyle{remark}
\newtheorem{remark}[thm]{Remark}
\newtheorem{fact}[thm]{Fact}
\newtheorem*{notation}{Notation}
\newtheorem*{greg}{$\Greg$}

\newcommand \seq[1]{{\langle{#1}\rangle}}
\newcommand{\uh}{\!\upharpoonright}
\newcommand{\da}{\!\downarrow\,}
\newcommand{\ua}{\uparrow}
\newcommand{\cat}{{}^\smallfrown}
\newcommand{\defn}[1]{{\emph{#1}}}

\newcommand{\Q}{\mathbb{Q}}
\newcommand{\dom}{\qopname\relax o{dom}}
\newcommand{\noqed}{\def\qedsymbol{}}

\renewcommand{\phi}{\varphi}
\renewcommand{\epsilon}{\varepsilon}
\renewcommand{\le}{\leqslant}
\renewcommand{\leq}{\leqslant}
\renewcommand{\nleq}{\nleqslant}
\renewcommand{\ge}{\geqslant}
\renewcommand{\geq}{\geqslant}

\newcommand{\s}{\sigma}
\newcommand{\e}{\epsilon}
\newcommand{\w}{\omega}
\newcommand{\cero}{\mathbf{0}}
\newcommand{\aaa}{\mathbf{a}}
\newcommand{\hh}{\mathbf{h}}
\newcommand{\D}{\mathcal{D}}
\newcommand \pp{\mathbf{p}}
\newcommand \qq{\mathbf{q}}

\newcommand \vphi{\varphi}

\newcommand \g{\gamma}
\newcommand \PP{\mathbb{P}}

\newcommand {\revmathfont}[1]{\mathsf{#1}}
\newcommand {\naughtsys}[1]{{\revmathfont{#1}}_0}
\newcommand\ACAz{\naughtsys{ACA}}
\newcommand\WKLz{\naughtsys{WKL}}
\newcommand\RCAz{\naughtsys{RCA}}
\newcommand\WWKLz{\naughtsys{WWKL}}
\newcommand\DNRz{\naughtsys{DNR}}
\newcommand {\Greg}{G_\delta\revmathfont{\text{-}REG}}

\newcommand {\Fin}{\text{\textit{Fin}}}

\begin{document}

%
%

\title{Uniform almost everywhere domination}

\author[Cholak]{Peter Cholak}

\address{Department of Mathematics\\ University of Notre Dame\\
  Notre Dame, IN 46556-5683}

\email{Peter.Cholak.1@nd.edu}

\author[Greenberg]{Noam Greenberg}

\address{Department of Mathematics\\ University of Notre Dame\\
  Notre Dame, IN 46556-5683}

\email{erlkoeing@nd.edu}

\author[Miller]{Joseph S.~Miller}

\address{Department of Mathematics \\
  University of Connecticut, U-3009 \\
  196 Auditorium Road\\
  Storrs, CT 06269}

\email{joseph.s.miller@gmail.com}

\thanks{Research partially supported by NSF Grants EMSW21-RTG-03-53748
  (all authors) and 02-45167 (Cholak). Miller's research was primarily
  supported by an NSF VIGRE Fellowship at Indiana University Bloomington.
  We wish to thank Bj{\o}rn
  Kjos-Hanssen and Denis Hirschfeldt for a number of interesting
  discussions about this work. We also thank Bj{\o}rn Kjos-Hanssen for
  allowing us to include Theorem~\ref{thm:Kjos-Hanssen} and its proof.}

\begin{abstract}
  We explore the interaction between Lebesgue measure and dominating
  functions. We show, via both a priority construction and a forcing
  construction, that there is a function of incomplete degree
  that dominates almost all degrees. This answers a question of Dobrinen
  and Simpson, who showed that such functions are related to the
  proof-theoretic strength of the regularity of Lebesgue measure for
  $G_\delta$ sets. Our constructions essentially settle the reverse
  mathematical classification of this principle.
\end{abstract}

\date{November 13, 2005}

\maketitle

\section{Introduction}

\subsection{Domination}

Fast growing functions have been investigated in mathematics for over
90 years. Set theorists, for example, have investigated the structure
$\w^\w/\Fin$ and the associated invariants of the continuum ever since
Hausdorff constructed his $(\w_1,\w_1^*)$-gap \cite{Hausdorff:gap};
today, this structure has a role to play in modern descriptive set
theory.

Fast growing functions have deep connections with computability. A
famous early example is that of Ackermann's function, defined in 1928
\cite{Ackermann}. This is a computable function that grows faster than
any primitive recursive function. This example was useful in
elucidating the mathematical concept of computability, an
understanding reflected in Church's Thesis.

In the 1960s, computability theorists became interested in functions
that grow faster than all computable functions.

\begin{definition} Let $f,g \colon \w\to \w$.
The function $f$ \emph{majorizes} $g$ if $f(n) \geq g(n)$ for all
$n\in\w$. If $f(n) \geq g(n)$ for all but finitely many $n$, then $f$
\emph{dominates} $g$. These are written as $f \geq g$ and $f \geq^*
g$, respectively. We call $f$ \emph{dominant} if it dominates all
(total) computable functions.
\end{definition}


Dominant functions were explored in conjunction with Post's Program.
The goal of Post's Program was to find a ``sparseness'' property of
the complement of a c.e.\ set $A$ that would ensure that $A$ is
incomplete. \citet{Yates:65} proved that even maximal c.e.\ sets,
which have the sparsest possible compliments among coninfinite c.e.\
sets, can be complete. This put an end to Post's Program, but not to
the study of sparseness properties.

Let $p_A(n)$ be the $n^\textit{th}$ element of the complement of $A$.
Having $p_A$ dominant would certainly imply that the complement of $A$
is sparse. On the other hand, \citet{Tennenbaum:63} and
\citet{Martin:66*1} showed that if $A$ is maximal, then $p_A$ is
dominant. Furthermore, Martin characterized the Turing degrees of both
the dominant functions and the maximal c.e.\ sets. He showed that there
is dominant function of degree $\aaa$ iff $\aaa$ is high (i.e., $\cero''
\leq \aaa'$), and that every high c.e.\ degree contains a maximal set.
Together, these results revealed a surprising connection between the
structure of c.e.\ sets, the place of their Turing degree within the
jump hierarchy, and domination properties of functions. Later research
explored further connections between domination properties, algebraic
properties and computational power.


In this paper, we consider the interaction between Lebesgue measure
and domination. Motivated by results on dominating functions in
generic extensions of set theory, \citet{MR2078930} introduced the
notion of a \emph{uniformly almost everywhere} (\emph{a.e.})
\emph{dominating} degree: a Turing degree $\aaa$ that computes a
function $f\colon\w\to\w$ such that
\[
\mu\;\left\{Z\in 2^\omega\,:\, (\forall g\in \omega^\omega)[\;g\leq_T
Z \implies g\leq ^* f\;]\right\}=1.
\]
(Here $\mu$ denotes the Lebesgue measure on $2^\w$.) We also call such
a function $f$ \emph{uniformly a.e. dominating}.

A natural goal is to characterize those Turing degrees that are
uniformly a.e.\ dominating. A function of degree $\cero'$ that dominates
almost all degrees was first constructed by
\citet[Theorem~4.3]{Kurtz:81}. (Kurtz used this result to exhibit a
difference between the $1$-generic and the (weakly) $2$-generic degrees:
the upward closure of the $1$-generic degrees has measure one
\cite[Theorem 4.1]{Kurtz:81}, while the upward closure of the (weakly)
$2$-generic degrees has measure zero \cite[Corollary 4.3a]{Kurtz:81}.)
Since the collection of uniformly a.e.\ dominating degrees is closed
upwards, Kurtz's result implies that every degree $\ge \cero'$ is in the
class. On the other hand, a uniformly a.e.\ dominating function is
dominant, and so by Martin's result, every uniformly a.e.\ dominating
degree is high. Thus, Dobrinen and Simpson asked whether either the
class of complete degrees (degrees above $\cero'$) or the class of high
degrees is identical to the class of uniformly a.e.\ dominating degrees.

Unfortunately, the truth lies somewhere in the middle. Binns,
Kjos-Hanssen, Lerman and Solomon \cite{Binns.KjosHanssen.ea:nd} showed
that not every high degree is uniformly a.e.\ dominating, or even
\emph{a.e.\ dominating}, an apparently weaker notion also introduced
by \citet{MR2078930}. They gave two proofs. First, by a direct
construction, they produced a high c.e.\ degree that is not a.e.\
dominating. (A similar result was independently obtained by Greenberg
and Miller, although their example was $\Delta^0_2$, not c.e.)

Second, \citet{Binns.KjosHanssen.ea:nd} showed that if $A$ has a.e.\
dominating degree, then every set that is $1$-random over $A$ is
$2$-random. If $A$ is also $\Delta^0_2$, then by \citet{Nies:nd*3},
$\emptyset'$ is $K$-trivial over $A$ and so $A$ is \emph{super-high}
(i.e., $A'\ge_{tt}\emptyset''$). By an index set calculation, there is
a c.e.\ set that is high but not \emph{super-high}, hence not a.e.\
dominating. It is open whether $\emptyset'$ being $K$-trivial over
$\mathbf{a}\leq \cero'$ implies that $\mathbf{a}$ is (uniformly) a.e.\
dominating; Kjos-Hanssen has some related results.

We prove that Dobrinen and Simpson's other suggested characterization
of the uniformly a.e.\ dominating degrees also fails.

\begin{thm}\label{thm: incomplete ce} There is an incomplete (c.e.)
  uniformly a.e.\ dominating degree.
\end{thm}

We provide two proofs of this result, although only one produces a
c.e.\ degree. In Section~\ref{priority} we use a priority argument to
construct an incomplete c.e.\ uniformly a.e.\ dominating degree and in
Section~\ref{forcing} we present a more flexible forcing construction
of an incomplete uniformly a.e.\ dominating degree.

\subsection{Domination and Reverse Mathematics}

As observed by \citet{MR2078930}, uniformly a.e.\ dominating degrees
play a role in determining the reverse mathematical strength of the
fact that the Lebesgue measure is regular.  For an introduction to
reverse mathematics, the reader is directed to \citet{Simpson:98}.

Regularity means that for every measurable set $P$ there is a
$G_\delta$ set $Q\supseteq P$ and an $F_\sigma$ set $S\subseteq P$
such that $\mu(S)=\mu(P)=\mu(Q)$, where a $G_\delta$ set is the
intersection of countably many open sets and an $F_\sigma$ set is the
union of countably many closed sets.  Hence the following principle is
implied by the regularity of the Lebesgue measure.

\begin{greg}
\em For every $G_\delta$ set $Q \subseteq 2^\omega$ there is an
$F_\sigma$ set $S \subseteq Q$ such that $\mu(S)=\mu(Q)$.
\end{greg}

Recall that the $G_\delta$ sets are exactly those that are $\Pi^0_2$
in a real parameter (that is, boldface $\pmb{\Pi}^0_2$), and the
$F_\sigma$ sets are exactly the $\pmb{\Sigma}^0_2$ sets. Hence we can
consider $\Greg$ as a statement of second order arithmetic. We will
see that $\Greg$, which appears to be a natural mathematical
statement, does not fall in line with the commonly occurring systems
of reverse mathematics. In particular, we examine the chain
\[
\RCAz \subsetneq \DNRz \subsetneq \WWKLz \subsetneq \WKLz \subsetneq
\ACAz.
\]
Here $\RCAz$ is the standard base system that all of the other systems
extend; $\WKLz$ is $\RCAz$ plus weak K\"{o}nig's lemma; and $\ACAz$ is
$\RCAz$ plus the scheme of arithmetic comprehension. These systems are
studied extensively in \cite{Simpson:98}. The system $\WWKLz$ is
somewhat less standard. It consists of $\RCAz$ plus ``weak weak
K\"{o}nig's lemma", which is introduced in \citet{YuSimpson:90}. A
large amount of basic measure theory can be proved in $\WWKLz$, so it
is a natural system for us to be concerned with. The final system,
$\DNRz$, is less natural from a proof-theoretic standpoint but very
natural for computability theorists. It is $\RCAz$ plus the existence
of a function that is diagonally non-recursive; see
\citet{GiustoSimpson} and \citet{Jockusch:NoFixedPoints}.

Kurtz's result that $\cero'$ is uniformly everywhere dominating
essentially shows that $\Greg$ follows from $\ACAz$. This relies on
the following:

\begin{thm}[Theorem~3.2 of \citet{MR2078930}]
  \label{thm: Dobrinen and Simpson}
  A Turing degree $\aaa$ is of uniformly a.e.\ dominating degree iff for every
  $\Pi^{0}_2$ set $Q \subseteq 2^\omega$ there is a $\Sigma^{0}_2(\aaa)$
  set $S \subseteq Q$ such that $\mu(S)=\mu(Q)$.
\end{thm}

Dobrinen and Simpson conjectured that $\Greg$ and $\ACAz$ are equivalent
over $\RCAz$ (\cite[Conjecture 3.1]{MR2078930}). This is not true; in
fact, there is an $\w$-model of $\Greg$ that omits $\cero'$ and hence is
not a model of $\ACAz$. This was discovered by B.\ Kjos-Hanssen after
the circulation of the priority-method proof of Theorem~\ref{thm:
incomplete ce}. This proof appears to be too rigid to allow us to obtain
a version with cone avoidance, but Kjos-Hanssen found a clever way to
build the $\w$-model without such a result. His construction is
presented in Section~\ref{Kjos-Hanssen}.

The forcing construction is flexible enough to prove cone avoidance
and more. We can thus improve Kjos-Hanssen's result by showing that
$\Greg$ does not imply even systems much weaker than $\ACAz$:

\begin{thm}
\label{thm: Greg is weak} $\RCAz+ \Greg$ does not imply $\DNRz$.
\end{thm}

But although $\Greg$ seems to lack proof-theoretic strength, none of
the traditional systems below $\ACAz$ are strong enough to prove it:

\begin{proposition}[Remark~3.5 of \citet{MR2078930}]
$\WKLz$ does not imply $\Greg$.
\end{proposition}

The proposition follows easily from the fact that there is an
$\w$-model of $\WKLz$ that consists of low sets; by formalizing
Theorem~\ref{thm: Dobrinen and Simpson}, every $\w$-model of $\Greg$
must include uniformly a.e.\ dominating degrees, which by Martin's
result are high.

Furthermore, $\Greg$ seems to be ``orthogonal" to the traditional
systems in that its strength is insufficient to lift one such system
to the system above it:

\begin{thm} \label{thm: Greg is orthgonal}
$\WKLz + \Greg$ does not imply $\ACAz$; $\WWKLz + \Greg$ does not
imply $\WKLz$.
\end{thm}

It remains open whether $\DNRz + \Greg$ implies $\WWKLz$.


\subsection{Notation, conventions and other technicalities}

Our computability theoretic notation is not always classical or
consistent, but hopefully completely understandable. Thus,
$\seq{\phi_e}_{e\in\omega}$ is an effective list of all Turing
functionals with oracle, and we write $\phi^f_e(x),
\phi^f_{e,s}(x)\da$, etc. This notation will be used when we try to
diagonalize against some oracle $f$ (so $\phi_e\colon \w^\w \to
\w^\w$). On the other hand, for domination purposes, we write Turing
functionals as $\Phi(Z;x)$ and $\Phi(Z;x)[s]$. In fact, we only need
to consider a single $\Phi$:

\begin{lemma}
  There is a partial computable functional $\Phi\colon 2^\omega\to
  \omega^\omega$ such that if
  \[
  \mu\;\{Z\in 2^\omega\,:\, \text{if $\Phi(Z)$ is total, then
  }\Phi(Z)\leq ^* f\}=1,
  \]
  then $f$ is uniformly a.e.\ dominating.
\end{lemma}

\begin{proof}
  Let $\seq{\Psi_i}_{i\in\omega}$ be an effective list of partial
  computable functionals $2^\omega\to \omega^\omega$ and define
  $\Phi(0^i1Z) = \Psi_i(Z)$.
\end{proof}

We assume that $\Phi$ has the following (standard) properties (for
every $s,n\in \w$ and $Z\in 2^\w$):
\begin{enumerate}
\item $\Phi(Z;n)\da[s]$ implies $\Phi(Z;n)[s]\leq s$.
\item $\Phi(Z;n)\da[s]$ implies $(\forall m<n)\;\Phi(Z;m)\da[s]$.
\end{enumerate}

We let $\dom \Phi$ be the collection of $Z$ such that $\Phi(Z)$ is
total. For $n\in \w$, we let $D_n = \{Z\in 2^\omega\,:\,
\Phi(Z;n)\da\}$. For a stage $s\in\omega$, $D_n[s]$ is given the
obvious meaning. For $g\in\omega^{\leq\omega}$, let
\[
D_{[n,m)}[g] = \left\{Z\in 2^\omega\,:\, (\forall k\in
  [n,m))\; \Phi(Z;k)\da[g(k)]\right\},
\]
(including the case where $m=\infty$). It follows from condition (1)
that if $Z\in D_{[n,m)}[g]$, then $g$ majorizes $\Phi(Z)$ on the
interval $[n,m)$.

\section{A proof of Theorem~\ref{thm: incomplete ce} via a priority
  construction} \label{priority}

In this section we prove Theorem~\ref{thm: incomplete ce}. We build
$f\colon\omega\to\omega$ by giving a computable sequence of
approximations $\seq{f_s}_{s\in\omega}$. Assuming the limit exists,
$f=\lim f_s$ is $\Delta^0_2$. To ensure that $f$ has c.e.\ degree, it is
enough to require that $f$ is approximated from below. Formally,
$(\forall n)(\forall s)\; f_s(n)\leq f_{s+1}(n)$. This means that $W =
\{\seq{n,m}\,:\, f(n)\ge m\}$ is a c.e.\ set; it is clear that
$f\equiv_T W$.

To ensure that $f$ is incomplete we will enumerate a c.e.\ set $B$ and
meet the requirement

$R_e\colon \phi^f_e\neq B$,

\noindent for each $e\in\omega$. These requirements will be handled by
incompleteness strategies. The same strategies are responsible for
assigning values to $f$, which essentially means that they must make
$f$ large enough to be uniformly a.e.\ dominating. This can be
accomplished if they are supplied with appropriate approximations to
the measure of $\dom \Phi$. These approximations are given by measure
guessing strategies.

We describe the incompleteness and measure guessing strategies first,
in relative isolation. Then we explain the priority tree and the full
construction.

\subsection{Incompleteness Strategy}

Let $\sigma$ be an agent assigned the goal of ensuring that
$\phi^f_e\neq B$, for some index $e=e(\sigma)$. When $\sigma$ is
initialized, it chooses a \emph{follower} $x=x(\sigma)$ that has not
been used before in the construction. A typical incompleteness
strategy would wait for a computation $\phi^f_e(x)\da = 0$, preserve
$f$ on the use of this computation, and enumerate $x$ into $B$. The
main difference is that our incompleteness strategy will be proactive:
it is permitted to change the values of $f$ to make $\phi^f_e(x)\da =
0$. Indeed, only the incompleteness agents change the values of $f$ at
all, so they are not only permitted to make these changes, it is
crucial that they do so.

Three restrictions are placed on $\sigma$'s ability to change the
values of $f$.  First, as already mentioned, it cannot decrease the
current values of $f$. Second, higher priority agents (who wish to
preserve diagonalizing computations) impose restraint $N=N(\sigma)$;
$\sigma$ is not allowed to change $f\uh N$. The third restriction
(which ensures that eventually $f$ will be dominating) involves a
rational parameter $\epsilon = \epsilon(\sigma)$. For $\sigma$ to
permanently protect a computation $\phi^f_e(x)\da = 0$ with use $r$,
it must be the case that
\begin{equation*}
  \label{eq:dfum}\tag{$\lozenge$} \mu\left(\dom \Phi
  \smallsetminus D_{[N,r)}[f]\right)\leq \epsilon.
\end{equation*}
In other words, $\sigma$'s action (in protecting $f\uh r$) prevents
$f$ from majorizing $\Phi(Z)$ above $N$ for no more than $\epsilon$ of
all $Z\in 2^\omega$.  This is the restriction that forces $\sigma$ to
increase the values of $f$.

The first two restrictions place no significant burden on $\sigma$,
but the third is more demanding. In fact, $\sigma$ cannot hope to meet
the third restriction without help because it does not know what
$\dom\Phi$ is. To approximate it, we supply $\sigma$ with two useful
pieces of information: a rational $q=q(\sigma)$ and a natural number
$M=M(\sigma)$ such that:
\begin{enumerate}
\item $q \leq \mu(\dom \Phi)$.
\item $\mu(D_M) \leq q + \epsilon/2$.
\end{enumerate}
In the full construction, these parameters are provided by a measure
guessing agent. If $\sigma$ is on the true path, then the values of
$q$ and $M$ that are supplied to $\sigma$ will meet conditions (1) and
(2).

We are now ready to describe the behavior of $\sigma$. The possible
states of $\sigma$ are {\sf active}, meaning that it is currently
imposing restraint to protect a computation $\phi^f_e(x)$, and {\sf
  passive}. When $\sigma$ is initialized, it is passive and it has
restraint $r(\sigma)=0$. If $\sigma$ ever becomes active, it will
remain so unless it is reset. This happens if the execution ever moves
left of $\sigma$, or if condition (2) proves to be false for either
$\sigma$ or a higher priority active agent. The details of the full
construction are below.

Say that $\sigma$ is visited at stage $s\in\omega$. If either $\sigma$
or a higher priority agent for $R_e$ is currently active, then there
is nothing to do.  Otherwise, $\sigma$ searches for a string $g\in
s^{<s}$ that has the following (computable) properties:
\begin{enumerate}
\item $g\supset f_s\uh N$;
\item $(\forall n\in [N,|g|))\; f_s(n) \le g(n)$;
\item $\mu(D_{[N,|g|)}[g]) > q - \e/2$; and
\item $\phi_{e,s}^g(x)\da =0$.
\end{enumerate}

If there is such a string $g$, then $\sigma$ lets $f_{s+1}\supset g$
and $r(\s) = |g|$. It enumerates $x$ into $B$ and declares itself {\sf
  active}. If there is no such $g$, then $\s$ does nothing and remains
{\sf passive}.

This completes the description of the incompleteness strategy. We
prove below that if $\sigma$ is on the true path and it ever becomes
satisfied, then (\ref{eq:dfum}) holds. Because agents that are not on
the true path might also attempt to protect computations, what we
actually prove is stronger: if an agent ever becomes {\sf active}
(hence is imposing restraint), either (\ref{eq:dfum}) holds or the
agent is eventually reinitialized (so that its restraint is removed).

\begin{remark}
  Unlike many tree constructions, it is important that at most one
  node on each level (i.e.\ at most one node per requirement) imposes
  restraint. Say a node at level $e$ ensures that $f$
  dominates except for a set of size at most $\e_e$. We will argue
  that $f$ dominates almost everywhere, using the fact that
  $\lim_{e\to \infty} \sum_{e'>e} \e_{e'} =0$. If several nodes on the
  same level $e$ were to impose restraint, then $\e_e$ must be counted
  more than once, making the calculation incorrect. This is why we
  stipulated that if $\sigma$ is visited at some stage $s$ and if at
  the same stage, some $\sigma'<_L \sigma$ on the same level is active,
  then $\sigma$ does not act. Of course, we are making use of the fact
  that $\sigma'$'s success is also $\sigma$'s.
\end{remark}

\subsection{Measure Guessing Strategy}

Measure guessing agents change neither $f$ nor $B$ and they impose no
restraint on other agents. Their only function is to provide the values
of $q$ and $M$ to the incompleteness agents at the next higher level. A
measure guessing agent $\tau$ is initialized with a rational parameter
$\delta=\delta(\tau)$. Its primary job is to find a rational $q$ that
approximates the measure of $\dom \Phi$ from below to within $\delta$.
This is done as follows. Divide the interval $[0,1]$ into subintervals
of length $\delta$. When $\tau$ is visited at stage $s$, it compares,
for each $n\leq s$, the measure of $D_n[s]$ with that of $D_n[t]$, where
$t$ was the previous stage at which $\tau$ was visited. If the measure
of some $D_n$ has crossed the threshold from one subinterval $I'$ to one
on its right $I$, then (for the least such $n$) $\tau$ guesses that $q =
\min I$ approximates the measure of $\dom \Phi$. Assume that $\tau$ is
visited infinitely often and $\min I$ is the largest approximation
guessed infinitely often. Then $\mu(D_n)\geq \min I$ for all $n\in\w$
and $\mu(D_n)>\max I$ for finitely many $n$. Therefore, $\mu(\dom \Phi)
\in I$.

We give the details. Let $d=\lceil 1/\delta\rceil$. The outcomes of
$\tau$ will be of the form $\seq{q,M}\in \Q\times\omega$, where
$q\in\{0,\delta,2\delta,\ldots,d\delta\}$. When $\tau$ is first
initialized, its outcome is $\seq{0,0}$. Say that $\tau$ is visited at
stage $s\in\omega$ and that the previous visit occurred at stage $t<s$.
To provide a guess, $\tau$ looks for $n\leq s$ and $b\leq d$ such that
$\mu(D_n[t])< b\delta$ but $\mu(D_n[s])\geq b\delta$. For the greatest
such $b$ (or equivalently, the $b$ corresponding to the least such $n$),
$\tau$ lets $q=b\delta$. Otherwise, $\tau$ lets $q=0$. Finally, $\tau$
takes the least $M$ such that $\mu(D_M[s]) < q+\delta$. Because
$\mu(D_n[s])$ is monotonically decreasing as a function of $n$, for all
$n\geq M$ we also have $\mu(D_n[s]) < q + \delta$. The outcome of $\tau$
at stage $s$ is $\seq{q,M}$.

\begin{remark}
  \label{rem:measure} Suppose that $\tau$ has outcome $\seq{q,M_0}$ at
  stage $s_0$ and outcome $\seq{q,M_1}$ at $s_1>s_0$. Further suppose
  that whenever $\tau$ is visited at a stage $t$ between $s_0$ and $s_1$,
  its outcome at $t$ is of the form $\seq{q',M'}$ with $q'\leq q$. Then
  $M_1 = M_0$.
\end{remark}

\subsection{The Priority Tree}

As usual, agents are organized on a tree, with the children of an
agent representing its potential outcomes. Write $\alpha\subset\beta$
to mean that $\beta$ is a proper extension of $\alpha$. Each agent
comes with a linear ordering $<_L$ on its children. We extend $<_L$ to
other nodes as follows: say that $\alpha$ is \emph{to the left} of
$\beta$ and write $\alpha<_L\beta$ if there are $\rho\subseteq\alpha$
and $\nu\subseteq\beta$ such that $\rho$ and $\nu$ have the same
parent and $\rho<_L\nu$. Write $\alpha < \beta$ if either
$\alpha\subset\beta$ or $\alpha<_L\beta$. This is the total ordering
lexicographically induced on the tree by the ordering we impose on the
children of agents. If $\alpha<\beta$, then we say that $\alpha$
\emph{has higher priority} than $\beta$.

The even levels of the priority tree are devoted to measure guessing
agents and the odd levels to incompleteness agents. A measure guessing
agent $\tau$ at level $2k$ is supplied with the parameter
$\delta(\tau)=3^{-k}/2$. As described above, its outcomes have the
form $\seq{q,M}\in \Q\times\omega$, where $q$ is restricted to
rationals of the form $b\delta(\tau)$. The outcomes are ordered first
by $q$ and then by $M$, with larger numbers \emph{to the left of}
smaller numbers.

An incompleteness agent $\sigma = \tau\cat\seq{q,M}$ at level $2k+1$ has
parameters $e(\sigma)=k$ and $\epsilon(\sigma)=3^{-k}=2\delta(\tau)$. We
also obviously set $q(\sigma)=q$ and $M(\sigma)=M$. The two final
parameters, the follower $x(\sigma)$ and the restraint $N(\sigma)$
imposed by stronger nodes, are determined when $\sigma$ is
\emph{initialized}. To initialize $\sigma$ at stage $s\in\omega$, set
its state to {\sf passive}, let the restraint $\s$ imposes $r(\sigma)=0$
and choose a follower $x(\sigma)\in\omega$ that has not yet been
assigned in the construction. Furthermore, set
\[
N(\sigma) = \max\{\,r(\sigma')\,:\, \sigma'<_L\sigma\text{ is active
  at stage $s$}\,\}.
\]
The children of $\sigma$ are
$\sigma\cat\textsf{active}<_L\sigma\cat\textsf{passive}$.

\subsection{Full Construction}

Let $f_0(n)=0$ for all $n\in\omega$. The construction proceeds in
stages. The preliminary phase of stage $s\in\omega$ involves
reevaluating, and possibly \emph{resetting}, currently active
incompleteness agents. Reset agents must be reinitialized the next
time they are visited. Say that $\sigma = \tau\cat\seq{q,M}$ is active
at stage $s$. If $\mu(D_M[s]) > q + \epsilon(\sigma)/2$, then $\sigma$
acted based on a false assumption and it could be the case that
$\sigma$ is forcing $f\uh r(\sigma)$ to remain prohibitively small.
Therefore, we reset $\sigma$. We also reset all previously initialized
incompleteness agents of lower priority than $\sigma$ (to allow them
to recompute their restraints the next time they are visited).

\begin{remark}
  Suppose that $\tau$ lies on the true path and that
  $\sigma = \tau\cat\seq{q,M}$ is active at stage $s$. Further suppose
  that $\tau$'s guess is found to be incorrect at $s$
  (in other words, $\mu(D_M[s]) > q + \delta(\tau)$). Then the next time that
  $\tau$ is accessible, its new outcome lies to the left of
  $\sigma$ and so $\sigma$ is reset. It would seem that this mechanism
  would suffice and that explicit resetting is unnecessary. However,
  unlike many tree constructions, we need to be concerned with the
  restraint imposed by nodes that lie to the left of the true
  path. Such unwarranted restraint may prevent $f$ from sufficiently
  dominating, and so needs to be reset when found incorrect.
\end{remark}

During the main phase of stage $s$, we execute the strategies of
finitely many agents on the priority tree, following a path of length
at most $s$. This is done in substages $t\leq s$. We begin at substage
$t=0$ by visiting the root node $\alpha_0 = \lambda$. Say that we are
visiting an agent $\alpha_t$ at substage $t$. First, reset any
incompleteness agents $\sigma$ such that $\alpha_t<_L\sigma$. (Note
that if $\sigma$ is reset and $\sigma<\sigma'$, then $\alpha_t <_L
\sigma'$, so $\sigma'$ is also reset.)

\emph{Case $1\colon$ $\alpha_t$ is a measure guessing agent.} If the
outcome of $\alpha_t$ at stage $s$ is $\seq{q,M}$, then let
$\alpha_{t+1} = \alpha_t\cat\seq{q,M}$ and end the substage.

\emph{Case $2\colon$ $\alpha_t$ is an incompleteness agent.} If
$\alpha_t$ has never been visited before or has been reset since the
last time it was visited, then it is initialized. If $\alpha_t$ is
currently active, then end the substage and set $\alpha_{t+1} =
\alpha_t\cat{\sf active}$. Similarly, if there is a higher priority
agent for $R_e$ that is active at stage $s$, then set $\alpha_{t+1} =
\alpha_t\cat{\sf passive}$ and end the substage. Otherwise, execute
the incompleteness strategy for $\alpha_t$ at stage $s$. If $\alpha_t$
becomes active (so that changes are made to $f$ and $B$), then end
stage $s$ entirely. Otherwise, let $\alpha_{t+1} = \alpha_t\cat{\sf
  passive}$ and end the substage.

This continues until substage $t=s$ is completed or until stage $s$ is
explicitly ended because an incompleteness agent becomes active.
Finally, for any $x<\dom f_s$, if not expressly altered by us during the
stage, we let $f_{s+1}(x)=f_s(x)$. This completes the construction.

\subsection{Verification}

Inductively define the \emph{true path} to be the leftmost path
visited infinitely often. In particular:
\begin{itemize}
\item The root node $\lambda$ is on the true path.
\item If $\rho$ is on the true path and $\nu$ is the leftmost child
  of $\rho$ that is visited infinitely often (if such exists), then
  $\nu$ is on the true path.
\end{itemize}
It is clear that if $\rho$ is on the true path, then there is a stage
$s\in\omega$ after which no agent left of $\rho$ is ever visited.

\begin{claim}
  \label{cl:true-inc} If $\sigma$ is an incompleteness agent on the
  true path, then there is a stage $s\in\omega$ at which $\sigma$ is
  initialized and after which it will never be reset.
\end{claim}
\begin{proof}
  Take a stage $t\in\omega$ large enough that no agent left of
  $\sigma$ will ever again be visited. By induction, we may also
  assume that $t$ is large enough that the agents
  $\sigma'\subset\sigma$ have all been initialized for the final time
  (and will never be reset). None of these $\sigma'$ can become active
  after stage $t$, or else the execution would move left of $\sigma$.

  Although no $\sigma'<_L\sigma$ can become active after stage $t$,
  they can be reset in the preliminary phase of the construction and
  this will reset $\sigma$. But only active agents become reset and
  only finitely many $\sigma'<_L\sigma$ are active at stage
  $t$. Therefore, there is a stage $t'\geq t$ after which no agents
  left of $\sigma$ are ever reset.

  This leaves only one way that $\sigma = \tau\cat\seq{q,M}$ can be
  reset at any stage $t''\geq t'$: if $\sigma$ is active at stage
  $t''$ and $\mu(D_M[t'']) > q + \epsilon(\sigma)/2$. But if this is
  the case, then $\seq{q,M}$ cannot be the outcome of $\tau$ after stage
  $t''$, contradicting the fact that $\sigma$ is on the true
  path. Therefore, $\sigma$ is never reset after stage $t'$. But
  $\sigma$ is visited infinitely often, so there is a stage
  $s\in\omega$ at which $\sigma$ is initialized and after which it
  will never be reset.
\end{proof}

\begin{claim} The true path is infinite.
\end{claim}
\begin{proof}
  We prove that there is no last node on the true path. First,
  consider an incompleteness agent $\sigma$ on the true path. By
  Claim~\ref{cl:true-inc}, there is a last stage $t$ at which $\sigma$
  is initialized. After stage $t$, $\sigma$ may become active at most
  once, so one of the outcomes of $\sigma$ is eventually
  permanent.

  Now consider a measure guessing agent $\tau$ on the true path. The
  first coordinate of the outcome of $\tau$ is taken from the finite
  set $Q = \{b\delta(\tau)\,:\, 0\leq b\leq \lceil
  1/\delta(\tau)\rceil\}$. Let $q$ be the greatest element of $Q$ that
  occurs as the first coordinate of the outcome infinitely
  often. Assume that no greater first coordinate occurs after stage
  $s\in\omega$. Let $\seq{q, M}$ be the outcome of $\tau$ at some stage
  $\geq s$. By Remark~\ref{rem:measure}, if $\seq{q',M'}$ is the outcome
  of $\tau$ at some other stage $\geq s$, then either $q'<q$ or $q'=q$
  and $M'=M$. Therefore, either $\seq{q,M} <_L \seq{q',M'}$ or
  $\seq{q,M}=\seq{q',M'}$, and the second case occurs infinitely
  often. Hence $\tau\cat\seq{q,M}$ is on the true path.
\end{proof}

\begin{remark}
  Let $\tau$ be a measure guessing node, and suppose
  that $\tau$ and $\tau\cat\seq{q,M}$ are on the true path.
  Then $\mu(D_n)\geq q$ for all $n\in\omega$ and
  $\mu(D_M) \leq q + \delta(\tau)$.
  Therefore, $\mu(\dom \Phi) \in [q,q+\delta(\tau)]$.
\end{remark}

We are primarily interested in the incompleteness agents that are
eventually permanently {\sf active}. Let the set of all such agents be
$G = \{\sigma_0,\sigma_1,\dots\}$, with $\sigma_0<\sigma_1<\cdots$.
The fact that we can thus enumerate $G$ relies on the following:

\begin{fact}
  The collection of nodes that lie either on, or to the left of the
  true path that are ever visited has order type $\w$ under $<$. This is
  because for each node $\alpha$ on the true path, only finitely many
  nodes to the left of $\alpha$ are ever visited.
\end{fact}

\begin{claim}
  \label{cl:G-active} Assume that $\sigma$ is initialized at stage
  $s\in\omega$ and is never reset after stage $s$. Suppose that
  $\sigma' < \sigma$. Then if $\sigma'$ is active at $s$, it
  remains permanently so (hence $\sigma'\in G$); otherwise, $\sigma'$
  never becomes active after $s$ (hence $\sigma'\notin G$).
\end{claim}

\begin{proof}
  First assume that $\sigma'$ is active at stage $s$. If $\sigma'$ is
  ever reset, then every lower priority agent is reset, including
  $\sigma$. But this never happens, so $\sigma'\in G$.

  Now suppose that $\sigma'$ is not active at stage $s$. It follows
  that $\sigma'\cat {\sf passive} < \sigma$ (as $\sigma$ is accessible
  at stage $s$). If $\sigma'$ becomes active at some later stage,
  then $\sigma' \cat{\sf active}$ would be accessible. But this
  would reset $\sigma$ because $\sigma'\cat{\sf active}$ lies to the
  left of $\sigma$.
\end{proof}

For all $i\in\omega$, let $N_i$ and $r_i$ denote the final values of
$N(\sigma_i)$ and $r(\sigma_i)$, respectively.

\begin{claim} For all $i\in\omega$:
  \label{cl:more-G}
  \begin{itemize}
  \item[(a)] Once $\sigma_i$ becomes permanently active, $f$ cannot
    change below $r_i$.
  \item[(b)] $N_{i+1} = r_i$.
  \end{itemize}
\end{claim}
\begin{proof}
  (a) Assume that $\sigma_i$ is permanently active after stage
  $s\in\omega$. From $s$ onwards, $\sigma_i$ imposes restraint $r_i$
  on weaker agents, so such agents do not change $f\uh r_i$. Any
  action by a stronger agent is impossible after the last stage $s_i$
  at which $\sigma_i$ is initialized, and $s_i<s$.

  (b) At stage $s_i$, the agents $\sigma < \sigma_i$ that are active
  are exactly $\sigma_0,\dots, \sigma_{i-1}$, and their restraints
  have reached their final values. Thus $\sigma_i$ defines $N_i = \max\{
  r_{j}\,:\, j<i\}$ at stage $s_i$. When $\sigma_i$ later becomes active, it
  imposes a permanent restraint $r_i$, which is greater than $N_i$.
  It follows that $r_0<r_1<\cdots$, and so $N_{i+1} = r_i$.
\end{proof}

\begin{claim}
  \label{cl:G-infinite} $G$ is infinite.
\end{claim}
\begin{proof}
  We can enumerate $\seq{\phi_e}_{e\in\w}$ in such a way that there are
  infinitely many $e$ such that for all $t\in \w$, for all $x\le t$
  and all $g\in {(t+1)}^{\le t}$ we have $\phi_{e,t}^g(x)\da=0$;
  we retroactively assume that we used such an enumeration. We
  will show that for each such $e$, $G$ contains an agent working for
  $R_e$.

  Pick such an $e$ and let $\sigma = \tau\cat\seq{q,M}$ be the agent of
  length $2e+1$ on the true path. Assume that the final initialization
  of $\sigma$ occurs at stage $s\in\omega$.

  \emph{Case $1\colon$ An agent $\sigma'<_L\sigma$ for $R_e$ is active
    at stage $s$.}  If $\sigma'$ is ever reset, then $\sigma$ would
  also be reset. This is impossible, so $\sigma'\in G$.

  \emph{Case $2\colon$ No such $\sigma'$ exists.} No
  $\sigma'<_L\sigma$ becomes active after stage $s$, so as
  long as $\sigma$ remains passive, its full strategy will be executed
  every time it is visited. At stage $s$, a follower $x$ is chosen and
  the final restraint $N$ is determined. By Claim~\ref{cl:more-G},
  $f\uh N$ is fixed after stage $s$.

  We know that $\seq{q,M}$ is the correct outcome of $\tau$, so
  $(\forall n)\; \mu(D_n) \geq q$ (recall that $\seq{D_n}$ is a decreasing sequence.)
   Let $v= \max\{N,M\}$. There is a
  $t_0$ such that $\mu(D_v[t_0]) > q - \e(\s)/2$. For any string $g\in
  \w^{v+1}$ extending $f\uh N$ such that $g(n)\ge t_0$ for all
  $n\in [N,v+1)$, we have $\mu(D_{[N,|g|)}[g]) > q - \e(\s)/2$.

  Consider a stage $t\geq \max\{t_0,x,v+1\}$ at which $\sigma$ is accessible. Let
  $g = f\uh N \cat \seq{t}^{v+1-N}$. Of course $f_t(n)\le t$ for all
  $n$, so by the assumptions on $e$, $\phi^g_{e,t}(x)\da=0$. Thus
  $g$ satisfies all the conditions that make it eligible to be picked
  as a new initial segment of $f$. It follows that if $\s$ did not act
  before stage $t$, then it does so and becomes permanently {\sf active}.
\end{proof}

\begin{claim}
  $f=\lim_s f_s$ exists.
\end{claim}
\begin{proof}
  Combining Claims~\ref{cl:more-G}(b) and \ref{cl:G-infinite}, the
  intervals $\{\,[N_i,r_i)\,\}_{i\in\omega}$ partition
  $\omega$. Furthermore, by Claim~\ref{cl:more-G}(a), $f$ is stable on
  $[0,r_i)$ once $\sigma_i$ becomes permanently active. Therefore,
  $\lim_s f_s(n)$ converges for all $n\in\omega$.
\end{proof}

\begin{claim}
  \label{cl:lozenge} For all $i$, $\mu\left(\dom\Phi\smallsetminus
    D_{[N_i,r_i)}[f]\right )\leq \epsilon(\sigma_i)$.
\end{claim}
\begin{proof}
  Assume for a contradiction that $\mu\left(\dom \Phi\smallsetminus
    D_{[N_i,r_i)}[f]\right)>\epsilon(\sigma_i)$. Let $\sigma_i =
  \tau\cat\seq{q,M}$. Take $s\in\omega$ to be the stage at which
  $\sigma_i$ becomes permanently active and let $g\in
  \omega^{<\omega}$ be the string that was used at that
  activation. So $r_i=|g|$ and $g\subset f$. This implies that
  $D_{[N_i,r_i)}[g] = D_{[N_i,r_i)}[f]$.  But of course,
  $\dom \Phi\subseteq D_M$. Therefore, $\mu\left(D_M\smallsetminus
  D_{[N_i,r_i)}[g]\right)>\epsilon(\sigma_i)$.

  By the definition of the incompleteness strategy,
  $\mu(D_{[N_i,r_i)}[g])>q-\epsilon(\sigma_i)/2$. Also $r_i>M$, so
  $D_{[N_i,r_i)}[g]\subseteq D_M$. Together with the conclusion of the
  previous paragraph, we have $\mu(D_M) > q+\epsilon(\sigma_i)/2$. But
  then $\mu(D_M[t])> q+\epsilon(\sigma_i)/2$, for any sufficiently
  large $t\in\omega$.  Therefore, $\sigma_i$ would be reset at the
  first phase of stage $t$, which is a contradiction.
\end{proof}

\begin{claim}
  \label{cl:u.a.e.} $f$ is uniformly a.e.\ dominating.
\end{claim}
\begin{proof}
  Fix $e\in\omega$. The construction ensures that at most one
  incompleteness agent at each level can be active at a time; hence at
  most one can belong to $G$. Thus there is an $i\in\omega$ large
  enough that $(\forall j\geq i)\; |\sigma_j|\geq 2e+1$.  Furthermore,
  $\sum_{j\geq i} \epsilon(\sigma_j) \leq \sum_{k\geq e} 3^{-k} =
  3^{-e+1}/2$ for this choice of $i$. By Claims~\ref{cl:more-G}(b) and
  \ref{cl:G-infinite}, the intervals $\{\,[N_j,r_j)\,\}_{j\geq i}$
  partition $[N_i,\infty)$. Therefore, if $Z\in\bigcap_{j\geq i}
  D_{[N_j,r_j)}$, then $f$ majorizes $\Phi(Z)$ above $N_i$. By
  Claim~\ref{cl:lozenge},
  \[
  \mu\bigg(\dom \Phi\smallsetminus \bigcap_{j\geq i}
  D_{[N_j,r_j)}[f]\bigg) \leq \sum_{j\geq i}
  \mu\left(\dom \Phi\smallsetminus D_{[N_j,r_j)}[f]\right)\leq
  \frac{3^{-e+1}}{2}.
  \]
  In other words, the set of $Z\in\dom \Phi$ such that $f$ fails to
  dominate $\Phi(Z)$ has measure at most $3^{-e+1}/2$. But
  $e\in\omega$ was arbitrary, so $f$ is uniformly a.e.\ dominating.
\end{proof}

\begin{claim}
  $f <_T \cero'$.
\end{claim}
\begin{proof}
  It is sufficient to prove that $B\nleq_T f$. Fix an index
  $e\in\omega$.

  \emph{Case $1\colon$ There is an $R_e$ agent $\sigma_i\in G$.} Let
  $s\in\omega$ be the last stage at which $\sigma_i$ becomes active
  and let $x_i = x(\sigma_i)[s]$. By Claim~\ref{cl:more-G}(a), this is
  done via $g = f_{s+1}\uh r_i = f\uh r_i$. Because $\sigma_i$ is
  activated, we know that $x_i\in B$ and $\phi_{e,s}^g(x_i) = 0$.
  Therefore, $\phi_e^f(x_i) = 0 \neq B(x_i)$.

  \emph{Case $2\colon$ There is no agent for $R_e$ in $G$.} Let
  $\sigma = \tau\cat\seq{q,M}$ be the incompleteness agent of length
  $2e+1$ on the true path.  Assume that $\sigma$ is initialized for
  the last time at stage $s\in\omega$. Let $x = x(\sigma)[s]$ and $N =
  N(\sigma)[s]$. Note that $x\notin B$, because $\sigma$ does not
  become active after stage $s$ (else $\sigma\in G$, so we would be in
  Case~1).  Assume, for a contradiction, that $\phi_e^f(x)= 0$. Take
  $g\in \omega^{<\omega}$ such that $|g| > \max\{M,N\}$, $g$ is an
  initial segment of $f$, and $\phi_e^g(x)= 0$.

  By Claim~\ref{cl:G-active}, $\sigma' < \sigma$ is active at stage
  $s$ iff $\sigma'\in G$. Choose $i\in\omega$ such that $\sigma_{i-1}
  < \sigma < \sigma_i$. In particular, $N=N_i$. Since $\s$ is on the
  true path, we have $\s\subset \s_j$ for all $j\geq i$. This shows
  that $(\forall j\geq i)\; |\sigma_j|>2e+1$. Now take $m\in\omega$
  large enough that $r_m\geq |g|$. Then, $\sum_{j\in [i,m]}
  \epsilon(\sigma_j) < \sum_{k> e} 3^{-k} = 3^{-e}/2$. By the same
  argument as given in Claim~\ref{cl:u.a.e.},
  $\mu\left(\dom \Phi\smallsetminus D_{[N,r_m)}[f]\right) < 3^{-e}/2$.
  Therefore, $\mu\left(\dom \Phi\smallsetminus D_{[N,|g|)}[g]\right) <
  3^{-e}/2$. We know that $\mu(\dom \Phi)\geq q$. This proves that
  $\mu(D_{[N,|g|)}[g]) > q-3^{-e}/2$.

  Let $t>s$ be a stage at which $\s$ is accessible that is large
  enough so that $\phi_{e,t}^g(x)\da =0$. There is nothing
  stopping $\s$ from acting at stage $t$, which is the desired
  contradiction.
\end{proof}

\section{Reverse mathematics I: avoiding cone avoidance}
\label{Kjos-Hanssen}

Although the above c.e.\ construction (Section~\ref{priority}) does
not seem to generalize to yield a cone avoidance result, Kjos-Hanssen
showed that it does have a reverse mathematical consequence.

\begin{thm}[Kjos-Hanssen]
  \label{thm:Kjos-Hanssen}
  There is an $\omega$-model of $\RCAz + \Greg$ that does not contain
  $\cero'$. Hence $\Greg$ does not imply $\ACAz$ over $\RCAz$.
\end{thm}

\begin{proof} We construct an ideal of Turing degrees that (as an
$\w$-model) satisfies $\Greg$ but does not contain $\cero'$. The ideal
is the downward closure of an increasing sequence $\aaa_1<  \hh_1 <
\aaa_2 < \hh_2 <\dots$. We let $\aaa_1 = \cero$ and let  $\hh_1$ be
the c.e.\ degree given by Theorem~\ref{thm: incomplete ce}. The degree
$\hh_1$ is high. In the structure  $\D [\hh_1,\hh_1' ]$ we can find
some $\aaa_2$ that is  $\textup{low}(\hh_1)$ and that joins $\cero'$
to $\hh_1' =  \cero''$ (\citet{PosnerRobinson}). Now in the structure
$\D [ \aaa_2,\aaa_2' ]= \D[\aaa_2,\cero'']$, a relativized version of
Theorem~\ref{thm: incomplete ce} yields a degree $\hh_2<\cero''$ that
is uniformly almost everywhere dominating over $\aaa_2$. We cannot
have $\hh_2 \ge \cero'$ because $\hh_2 \ge \aaa_2$ and $\hh_2 <
\cero''$.

We now repeat. Again, using a relativized version of
\cite{PosnerRobinson}, we get an $\aaa_3 \in \D[\hh_2,\cero''']$ that
is $\textup{low}(\hh_2)$ and joins $\cero''$ to $\cero'''$; and an
$\hh_3\in \D[\aaa_3,\cero''')$ that is uniformly a.e.\ dominating over
$\aaa_3$. As before, $\hh_3$ is not above $\cero''$. But as $\hh_3 \ge
\aaa_2$ and $\cero' \vee \aaa_2 = \cero''$ we cannot have $\hh_3 \ge
\cero'$. The process now repeats itself to get the rest of the
sequence.
\end{proof}

\section{A proof of Theorem~\ref{thm: incomplete ce}
via a forcing construction}\label{forcing}

In this section we introduce a forcing notion that produces a
uniformly a.e.\ dominating function and that allows us to obtain cone
avoidance and more.

\subsection{The notion of forcing}

We approximate a function $f^G$. A \defn{condition} is a pair
$\seq{f,\e}$ where $f\in \w^{<\w}$ and $\e$ is a positive rational.
The idea is that $\pp=\seq{f,\e}$ states that $f$ is an initial
segment of $f^G$ and further $\pp$ makes an \defn{$\e$-promise}: the
collection of $Z\in \dom \Phi$ such that $f^G$ fails to majorize
$\Phi(Z)$ from $|f|$ onwards has size $<\e$. Thus, an extension
$g\supset f$ \defn{respects the $\e$-promise} if
\[ \mu \left(\dom \Phi \smallsetminus D_{[|f|,|g|)}[g] \right) < \e. \]

However, this is not a good definition of a partial ordering on the
conditions; we can have $g$ keep the $\e$-promise of $\seq{f,\e}$ and
$h$ keep the $\delta$-promise of $\seq{g,\delta}$ but fail to respect
the $\e$-promise of $\seq{f,\e}$. Thus, the relation would not be
transitive.
A simple modification ensures that every $h$ that keeps the
$\delta$-promise of $\seq{g,\delta}$ also keeps the $\e$-promise of
$\seq{f,\e}$. We say that a condition $\seq{g,\delta}$
\defn{extends} another condition $\seq{f,\e}$ if $f\subset g$, $\delta\le
\e$ and further, if $f\ne g$, then
\[ \mu\left(\dom \Phi\smallsetminus D_{[|f|,|g|)}[g] \right) + \delta <
\e.\]

\begin{lemma} The extension relation is transitive. \end{lemma}

\begin{proof}
  Suppose that $\seq{g,\delta}$ extends $\seq{f,\e}$ and is extended by
  $\seq{h,\g}$; we show that $\seq{h,\g}$ extends $\seq{f,\e}$. If either $f=g$ or
  $g=h$, then this is easy. Otherwise, the point is that
  \[ D_{[|f|,|h|)}[h] = D_{[|f|,|g|)}[g]\cap D_{[|g|,|h|)}[h]\]
  and so
  \begin{multline*}
    \mu\left(\dom\Phi \smallsetminus D_{[|f|,|h|)}[h]\right)\le
    \mu\left(\dom\Phi \smallsetminus D_{[|f|,|g|)}[g]\right) +
    \mu\left(\dom \Phi\smallsetminus D_{[|g|,|h|)}[h]\right)\le \\
    (\e-\delta) + (\delta- \g) = \e - \g,
  \end{multline*}
  as required.
\end{proof}

\begin{notation} We let $\PP$ be the collection of all conditions.
  For a condition $\pp=\seq{f,\e}$ we write $f^\pp = f$ and $\e^\pp =
  \e$. We also let $n^\pp = |f^\pp|$.
\end{notation}

\begin{lemma} For all $n<\w$, the set $ \{\pp\in \PP\,:\, n^\pp > n\}
  $ is dense in $\PP$. \label{lem: one}
\end{lemma}
\begin{proof}
  Let $\pp\in \PP$. Let $n> n^\pp$. For large enough $s$,
  \[\mu(D_n \smallsetminus D_n[s]) < \e^\pp.\] Now take
  $g\in \w^n$ extending $f$ such that $D_n[s] \subset
  D_{[n^\pp,n)}[g]$ (for example, by defining $g(m)=s$ for $m\ge n^\pp$).
  As $\dom\Phi \subset D_n$, we get that $\mu\left(\dom \Phi
  \smallsetminus D_{[n^\pp,n)}[g]\right) < \e^\pp$. We can then pick some small
  $\delta$ so that $\seq{g,\delta}$ extends $\pp$.
\end{proof}

If $G\subset \PP$ is generic (from now, by the word ``generic" we mean,
``sufficiently generic for the given argument"), then we let
\[f^G = \bigcup_{\pp\in G} f^\pp.\] The following is a corollary of
Lemma~\ref{lem: one}:

\begin{corollary}
  If $G$ is generic, then $f^G \in \w^\w$.
\end{corollary}

We now show that the $\e$-promises are kept.

\begin{lemma} Let $\pp\in \PP$, and suppose that $\pp\in G$ and that
  $G$ is generic.  Then
  \[ \mu \left(\dom \Phi \smallsetminus
    D_{[n^\pp,\w)}\left[f^G\right]\right) \le \e^\pp. \]
\end{lemma}

\begin{proof}
The sequence $\seq{D_{[n^\pp,m)}[f^G]}_{m> n^\pp}$ decreases with $m$
and
\[
D_{[n^\pp,\w)} \left[f^G\right] = \bigcap_{m>n^\pp}
D_{[n^\pp,m)}\left[f^G\right].
\]
So it is enough to prove that $\mu \left(\dom \Phi \smallsetminus
D_{[n^\pp,m)}\left[f^G\right]\right) \le \e^\pp$, for all $m>n^\pp$.
For any $m$, there is a $\qq\in G$ extending $\pp$ such that $n^\qq
\geq m$. By the definition of our partial ordering,
\[
\mu\left(\dom\Phi\smallsetminus
D_{[n^\pp,n^\qq)}\left[f^\qq\right]\right) < \e^\pp.
\]
But $D_{[n^\pp,n^\qq)}\left[f^\qq\right] \subseteq
D_{[n^\pp,m)}\left[f^G\right]$ because $f^\qq\subset f^G$, which
completes the proof.
\end{proof}

The following is immediate.

\begin{lemma} For all $\e>0$, the set $ \{\pp\in \PP\,:\, \e^\pp <
  \e\} $ is dense in $\PP$. \qed \label{lem: two}
\end{lemma}

As a corollary,

\begin{corollary} If $G\subset \PP$ is generic, then $f^G$ is uniformly
  almost everywhere dominating.
\end{corollary}

\subsection{Cone avoidance, etc.}

We show that if $G$ is generic, then indeed $f^G$ has no special
properties beyond domination. The following is the crucial technical
lemma. Consider the proof that if $g$ is Cohen generic over $A$ and
$A$ is not computable, then $g$ does not compute $A$. If some
condition $\tau\in 2^{<\w}$ forces that $\phi^g_e=A$ (and in
particular is total), then $A = \bigcup_{\s\supset\tau}\phi^\s_e $ is
computable because the collection of extensions of $\tau$ is
computable. We would like to do the same, but our partial ordering is
not computable. This difficulty is overcome as follows: given $\pp\in
\PP$, we can make a promise $\e^*$ much tighter than $\e^\pp$ and find
a rational $q$ sufficiently close to $\dom \Phi$ such that every
sufficiently long string $g\supset f^\pp$ respecting the
$\e^*$-promise satisfies $\mu(
  D_{[n^\pp,|g|)}[g])>q$ and every string satisfying the latter
(computable) condition respects the $\e^\pp$-promise. We can now
imitate the diagonalization argument (and more): if $\pp$ forces that
$\phi^{f^G}_e = A$, then we compute $A$ by examining $\phi^{g}_e$ for
strings $g$ satisfying the middle condition above. We argue that this
must give us all of $A$, for otherwise we could extended $\pp$ to keep
the $\e^*$-promise and avoid $\phi^{f^G}_e = A$.

\begin{lemma} Let $\pp\in \PP$. Then there is a c.e.\ set
  \[S \subset \{ f^\qq\,:\, \qq\le \pp\}\] and a $\pp^* \le
  \pp$ such that $\{ \qq\le \pp^*\,:\, f^\qq\in S\}$ is dense below
  $\pp^*$. \label{lem: tech}
\end{lemma}

\begin{proof}
  Find some $n>n^\pp$ such that $\mu\left(D_n \smallsetminus \dom
    \Phi\right) < \e^\pp/ 2$; also find a rational
  $q<\mu(\dom\Phi)$ such that $\mu(\dom \Phi)-q < \e^\pp/2$. Let
  \[ S = \left\{ g\in \w^{<\w}\,:\, \text{$g\supset f^\pp$, $|g|> n$,
  and $\mu(D_{[n^\pp,|g|)}[g]) > q$}\right\} \]

  It is clear that $S$ is c.e. Let $g\in S$; we
  show that for some $\qq\le \pp$ we have $f^\qq = g$. We have
  $\mu(D_{[n^\pp,|g|)}[g]) > \mu(\dom \Phi) - \e^\pp/2$ and
  $\mu(D_{|g|}) < \mu(\dom \Phi) + \e^\pp/2$; together we
  get $\mu\left(D_{|g|}\smallsetminus D_{[n^\pp,|g|)}[g]\right) <
  \e^\pp$. Of course, $\dom \Phi \subset D_{|g|}$ and so $\mu\left(\dom
    \Phi \smallsetminus D_{[n^\pp,|g|)}[g]\right) < \e^\pp$.

  Next, let $\pp^* = \seq{f^\pp, \delta}$ where $\delta < \e^\pp$ (so
  $\pp^*\le \pp$) and $\delta< \mu(\dom \Phi)-q$. Suppose that $\qq\le
  \pp^*$ and $n^\qq > n$. Then from
  \[ \mu \left( \dom\Phi \smallsetminus D_{[n^\pp,n^\qq)}\left[f^\qq\right]
  \right) < \delta \] we can conclude that
  $\mu(D_{[n^\pp,n^\qq)}[f^\qq]) > q$, so $f^\qq\in S$.
\end{proof}

\begin{lemma}
  If $A$ is noncomputable and $G\subset \PP$ is generic over $A$, then $f^G\not\ge_T A$.
\end{lemma}

\begin{proof}
  Let $\Psi\colon \w^\w\to 2^\w$ be a Turing functional. We show that
  the union of
  \begin{align*}
  E_0 & = \{ \pp\in \PP\,:\, \Psi(f^\pp) \perp A\}\quad\text{and}\\
  E_1 & = \{ \pp\in \PP\,:\, (\exists x)(\forall \qq\le
  \pp)\; \Psi(f^\qq,x)\ua\,\}
  \end{align*}
  is dense in $\PP$. Of course if
  $G\cap (E_0 \cup E_1) \ne \emptyset$, then $\Psi(f^G)\ne A$.

  Let $\pp\in \PP$, and take $S$ and $\pp^*$ given by
  Lemma~\ref{lem: tech}. If there are $g,g'\in S$ such that $\Psi(g) \perp \Psi(g')$,
  then one of them is incompatible with $A$, so $\pp$ has an extension
  in $E_0$. If $\bigcup_{g\in S} \Psi(g)$ is total, then it
  is computable, hence different from $A$. Again, $\pp$ has an
  extension in $E_0$.

  Otherwise, for some $x\in\w$, we have $\Psi(g,x)\ua$ for all $g\in
  S$. This implies that $\pp^*\in E_1$: for all $\qq\le \pp^*$,
  $f^\qq$ has an extension in $S$, and so $\Psi(f^\qq,x)\ua$.
\end{proof}

\begin{lemma}
  If $G\subset \PP$ is generic, then $f^G$ does not have PA-degree.
\end{lemma}

\begin{proof}
  Let $\psi\colon \w\to 2$ be a partial computable function that has
  no total computable extension. We show that $f^G$ does not compute a
  0-1 valued total extension of $\psi$.

  Let $\Theta\colon \w^\w\to 2^\w$ be a Turing functional. We show
  that the union of
  \begin{align*}
  E_0 & = \{ \pp\in \PP\,:\, (\exists x\in \dom \psi)\; \Theta(f^\pp,x)\da
  \ne \psi(x) \}\quad\text{and}\\
  E_1 & = \{ \pp\in \PP\,:\, (\exists x)(\forall \qq\le
  \pp)\; \Theta(f^\qq,x)\ua\,\}
  \end{align*}
  is dense in $\PP$. Of course if
  $G\cap (E_0 \cup E_1) \ne \emptyset$, then $\Theta(f^G)$ is not a total
  extension of $\psi$.

  Let $\pp\in \PP$; take $S$ and $\pp^*$ given by Lemma~\ref{lem:
    tech}. If there is a $g\in S$ such that $\Theta(g)\perp \psi$,
  then $\pp$ has an extension in $E_0$. If for some $x$, $\Theta(g,x)\ua$
  for all $g\in S$, then $\pp^*\in E_1$.

  One of the above must be the case; otherwise, we could compute a
  completion of $\psi$ as follows: for each $x$,
  search for a $g\in S$ such that $\Theta(g,x)\da$. For the first such $g$
  found, let $h(x) = \Theta(g,x)$. Then $h$ is computable, and must extend
  $\psi$.
\end{proof}

In fact, the same proof gives us somewhat more:

\begin{lemma}
  If $G\subset \PP$ is generic, then $f^G$ does not have DNR-degree.
\end{lemma}

\begin{proof}
  Let $\seq{\vphi_e}_{e\in\omega}$ be an enumeration of all partial computable
  functions from $\w$ to $\w$.

  Let $\Psi\colon \w^\w \to \w^\w$ be a Turing functional. We show
  that the union of
  \begin{align*}
  E_0 & = \{ \pp\in \PP\,:\, (\exists e)\; \Psi(f^\pp,e)\da =
  \vphi_e(e)\da\}\quad\text{and}\\
  E_1 & = \{ \pp\in \PP\,:\, (\exists x)(\forall \qq\le
  \pp)\; \Psi(f^\qq,x)\ua\,\}
  \end{align*}
  is dense in $\PP$. Of course if
  $G\cap (E_0 \cup E_1) \ne \emptyset$, then $\Psi(f^G)$ is not DNR.

  Let $\pp\in \PP$, and take $S$ and $\pp^*$ given by Lemma~\ref{lem:
    tech}. If there is a $g\in S$ and $e\in\omega$ such that
  $\Psi(g,e)\da = \vphi_e(e)\da$, then $\pp$ has an extension in $E_0$. If
  there is an $x$ such that $\Psi(g,x)\ua$ for all $g\in S$, then
  $\pp^*\in E_1$.

  Otherwise, define a total function $h\colon \w\to \w$ as follows:
  for each $x$, search for a $g\in S$ such that $\Psi(g,x)\da$.
  Let $h(x) = \Psi(g,x)$ for the first such $g$ discovered.
  Then $h$ is computable and DNR, which is impossible.
\end{proof}

\subsection{Relativization}

Let $B\subset \w$. All the results of this section relativize to
working above $B$.  Namely, we can define a notion of forcing $\PP_B$;
all is exactly as above, except that instead of $D_n$ and $\dom \Phi$
we use $\{Z\,:\, \Phi(B\oplus Z,n)\da\}$ and $\{ Z\,:\, \Phi(B\oplus
Z)\text{ is total}\}$. With exactly the same proofs, we see that a
generic yields a function $f^G$ that is uniformly almost everywhere
dominating over $B$. Lemma~\ref{lem: tech} now becomes the following:

\begin{lemma} Let $\pp\in \PP_B$. Then there is a set $S$, c.e.\ in
  $B$, such that $S\subset \{f^\qq\,:\, \qq\le \pp\}$
  and some $\pp^*\le \pp$ such that $\{\qq\le \pp^*\,:\, f^\qq\in S\}$
  is dense below $\pp^*$.
\end{lemma}

These are the analogous corollaries:
\begin{lemma} Suppose that $A\not\le_T B$ and that $G\subset \PP_B$ is
  generic over $A$. Then $B\oplus f^G \not\ge_T A$. \label{lem: rel
    cone avoision}
\end{lemma}

\begin{lemma} Suppose that $B$ does not have PA-degree and that
  $G\subset \PP_A$ is generic over $B$. Then $B\oplus f^G$ does not
  have PA-degree. \label{lem: rel PA}
\end{lemma}

\begin{lemma} Suppose that $B$ is not DNR, and that $G\subset \PP_B$
  is generic. Then $B\oplus f^G$ is not DNR.\label{lem:DNR}
\end{lemma}

\section{Reverse mathematics II}
\label{reverse}

The above forcing argument directly yields the results concerning the
proof-theoretic strength of $\Greg$.

Recall from \citet{Simpson:98} that $M\subseteq 2^\w$ is an $\w$-model
of $\RCAz$ iff it forms an ideal in the Turing degrees, and it is an
$\w$-model of $\WKLz$ iff it is a \emph{Scott system}: i.e., a Turing
ideal such that for all $A\in M$, there is a $B\in M$ of PA-degree
relative to $A$. Similarly, \citet{YuSimpson:90} proved that a Turing
ideal $M\subseteq 2^\w$ is an $\w$-model of $\WWKLz$ iff for all $A\in
M$, there is a $B\in M$ that is sufficiently random over $A$ (it is
enough that $B$ is $1$-random relative to $A$ by a result of
\citet{Kucera:85}).

\begin{proof}[Proof of Theorem~\ref{thm: Greg is weak}]
  An ideal of Turing degrees that models $\Greg$ but does not include
  any DNR degrees is easily built using Lemma~\ref{lem:DNR}.
\end{proof}

\begin{proof}[Proof of Theorem~\ref{thm: Greg is orthgonal}]
For the first part, we can inductively construct an $\w$-model of
$\WKLz + \Greg$ that avoids $\cero'$ by alternatively appealing to
Lemma~\ref{lem: rel cone avoision} and to the fact that a similar cone
avoidance lemma holds for obtaining paths through trees, hence for
PA-degrees (\citet[Theorem 2.5]{Jockusch.Soare:72*1}).

For the second part, a similar construction yields an $\w$-model of
$\WWKLz + \Greg$ that does not satisfy $\WKLz$, this time using
Lemma~\ref{lem: rel PA} and the following claim, which essentially
appears in \citet{YuSimpson:90}.

\begin{claim}
Suppose that $B$ does not have PA-degree. If $A$ is sufficiently
random over $B$, then $A\oplus B$ does not have PA-degree.
\end{claim}

\begin{proof}
This is from \cite[Page 172]{YuSimpson:90}. Let $E$ and $F$ be
disjoint c.e.\ sets that cannot be separated by any set computable in
$B$. By relativizing a result from
\citet{JockuschSoare:DegreesOfTheories}, the measure of
\[
S = \{Z\,\colon\, (\exists Y\leq_T Z\oplus B)\; E\subseteq Y\wedge
F\cap Y = \emptyset\}
\]
is zero. This is the collection of sets $Z$ such that $Z\oplus B$
computes a separator of $E$ and $F$. If $A$ is sufficiently random
over $B$, then $A\notin S$, meaning that it satisfies the claim. (In
fact, since $S$ is a $\Sigma^0_3(B)$-class, it suffices for $A$ to be
(weakly) $2$-random relative to $B$ \cite{Kurtz:81}.)
\end{proof}
\noqed\end{proof}




\bibliographystyle{cholak}
\bibliography{aedominating}

\end{document}